\documentclass[a4paper]{extarticle}
\usepackage[utf8]{inputenc}
\usepackage[T1]{fontenc}
\usepackage[backend=biber,minnames=2]{biblatex}
\usepackage{parskip} %or set \parindent to 0pt
\usepackage{graphicx}
\usepackage{fullpage}
\usepackage{color}
\usepackage{amsmath}
\usepackage{amssymb}
\usepackage{siunitx}
\usepackage{comment}
\usepackage[section]{placeins}
\usepackage{tikz}
\usepackage{arrayjobx}
\usepackage{hyperref}

\begin{document}
\title{Calculation of Iterated-Integral Signatures and Log Signatures}
\author{Jeremy Reizenstein, Centre for Complexity Science\thanks{Supported by the Engineering and Physical Sciences Research Council}\
%\with\\Ben Graham, Department of Statistics and Centre for Complexity Science
\\University of Warwick}
%\keywords{Free Lie algebra, Baker-Campbell-Hausdorff, Log signature}
\date{May 2015, revision of December 2017}
\maketitle

\begin{abstract}
We explain the algebra needed to make sense of the log signature of a path, with plenty of examples. We show how the log signature can be calculated numerically, and explain some software tools which demonstrate it.
\end{abstract}
\section{Introduction}

An iterated integral signature (or Brownian signature) is a sequence of numbers derived from a continuous path. The signature is important in the theory of differential equations driven by rough paths, and is finding applications in machine learning \cite{OxSigIntro}.

When trying to classify types of paths with a machine learning algorithm, if each path is specified as a long list of points then it can be hard for paths which are similar to be seen as similar by the algorithm. If instead the space is discretized and the path is specified in terms of the points hit, then the representation will be very large and sparse, which will make learning slow.
The signature may be a good way to represent a path to avoid these problems.
Similar smooth curves have similar signatures.
The signature is invariant on translations of the path. It is also invariant on reparametrizing the path. 

If two paths are different then their signatures will be different, unless one contains a section where the path exactly retraces itself.\cite{HL} A signature cannot distinguish between two paths where the only difference is that one path contains an extra section which is backtracked over. In many situations, the relevant paths will not contain precise backtracking.
It is possible that approximate backtracking might lead to numerically unstable signatures. If there is a dimension that is monotone along each path (for example time) then backtracking is impossible.

The log signature is a compressed form of the signature, which is more complicated to define. For the same amount (number of levels) of the signature, the log signature contains the same information in fewer numbers. The link between the two is continuous.

This note is a self-contained explanation of the operations required to calculate the log signature of a piecewise-linear path, and an explanation of some python code which assists in so doing.
Statements are made without proofs, and there is an emphasis on small examples and some explicit calculations. The main reference on the algebra is \cite{FLA}, and the application to paths is explained in chapter~2 of \cite{Lyons07}. The signature and log signature are formally defined in Sections~\ref{sec:sig} and \ref{sec:logsig}. Section~\ref{sec:code} explains the code and how to use it. Section~\ref{sec:vis} describes a visualisation of the log signature which may be helpful.

\section{Words}
\def \alph#1{{\color{blue}\mathbf{#1}}}
\def \lex{<_L}
Consider a set $\Sigma=\{\alph{1},\alph{2},\dots,\alph{d}\}$ of $d$ "letters" which has an ordering $<$. %, such as the set of upper case letters $A,\dots,Z$ or the numbers $\alph{1},\alph{2},\alph{3}$. 
The set of words with entries in $\Sigma$ is called the \emph{Kleene Star} of $\Sigma$ and is denoted by $\Sigma^*$. The length of a word $u$ is denoted $|u|$, for example $|\alph{3231}|=4$. The empty word is denoted by $\epsilon$ and the concatenation of words $u$ and $v$ is written $uv$%; concatenation is an associative operation
. If a word $w$ is equal to $uv$ for some words $u$ and $v$, then $u$ is said to be a \emph{prefix} of $w$ and $v$ is said to be a \emph{suffix} of $w$. If $u$ and $v$ are both not empty then they are said to be a \emph{proper} prefix and suffix of $w$. For example $\alph{1}$ is a proper suffix of $\alph{3231}$, and a suffix but not a proper suffix of $\alph{1}$.
The ordering $<$ on $\Sigma$ can be extended to an ordering $\lex$ on $\Sigma^*$ called \emph{alphabetical order} or \emph{lexicographic order} in the usual way. (Specifically: $\epsilon\lex u$ if $|u|>0$. For letters $a$ and $b$ and words $u$ and $v$, $au\lex bv$ if $a<b$ or both $a=b$ and $u\lex v$.)

%\section{Free vector space}
The \emph{free (real) vector space} on the finite set $\Sigma$ is the real vector space with basis given by the elements of $\Sigma$. We will just call it $\mathbb{R}^d$. An element looks like $a_1\alph{1}+\dots+a_d\alph{d}$ for real numbers $a_1,\dots,a_d$.

%\section{Tensor algebra}
The \emph{tensor algebra} of the vector space $\mathbb{R}^d$,  $T(\mathbb{R}^d)$, is the set of functions from $\Sigma^*$ to $\mathbb{R}$, or equivalently the set of (possibly) infinite formal sums of real multiples of words.  The word $u$ in $\Sigma^*$ is identified with the function which takes $u$ to 1 and all other words to 0, or the expression $1u$. 
We are only ever interested in finite restrictions of this in order to do calculations, in particular we choose an integer $m$ and ignore all words with length longer than $m$. $T^{\underline m}(\mathbb{R}^d)$ is the real vector space with basis given by words of length $m$ or less\footnote{this is known as $T^{(m)}(\mathbb{R}^d)$ in the notation of \cite{Lyons07}}. The function on words which returns their concatenation if it has length $m$ or less and returns the zero element otherwise extends uniquely to a bilinear operation\footnote{i.e. linear in each argument} on $T^{\underline m}(\mathbb{R}^d)$, which is called the \emph{concatenation product}. For example, in $T^{\underline 4}(\mathbb{R}^3)$,
\[(7\,\alph{132}+9\epsilon)(2\,\alph{1}+4\,\alph{21})=14\,\alph{1321}+18\,\alph{1}+36\,\alph{21}.\]

The exponential series defines a function from those elements of $T^{\underline m}(\mathbb{R}^d)$ which have no term in $\epsilon$ to $T^{\underline m}(\mathbb{R}^d)$: \begin{equation} \exp(x):=1\epsilon+x+\frac12x^2+\dots+\frac1{m!}x^m = \sum_{i=0}^m\frac1{i!}x^i
\end{equation}
where powers denote the concatenation product of $x$ with itself. Stopping the sum at $m$ is the same as stopping later than $m$, because all higher powers are $0$. $\exp(x)$ always has 1 as its $\epsilon$ component.

The series for $\log(1+x)$ defines a function $\log$ from those elements of $T^{\underline m}(\mathbb{R}^d)$ which have $\epsilon$ component 1 to $T^{\underline m}(\mathbb{R}^d)$:\begin{equation}
\log(1\epsilon+x):=x-\frac12x^2+\frac13x^3-\dots+(-1)^{m+1}\frac1mx^m
\end{equation}
There is never an $\epsilon$ component in the value of $\log(x)$. In fact, $\log$ and $\exp$ are inverses. For example, up to level 2,
\begin{align*} \exp(\alph{1}+3\alph{2})&=1\epsilon+\alph{1}+3\alph{2}+\frac12\alph{11}+\frac32(\alph{12}+\alph{21})+\frac92\alph{22}\\
\log(\exp(\alph{1}+3\alph{2}))&=\alph{1}+3\alph{2}.
\end{align*} 
%And up to level 3
%\begin{align*} \exp(\alph{1}+\alph{12})&=1\epsilon+\alph{1}+\frac12\alph{11}+\alph{12}
%+\frac16\alph{111}+\frac12(\alph{112}+\alph{121})\\
%\end{align*} 
%\newpage

\subsection{Lyndon words}
\def\necklace#1#2{
\newarray\Labels
\readarray{Labels}{#1}
\begin{tikzpicture}
%code based on http://www.texample.net/tikz/examples/cycle/
\def \n {#2}
\def \radius {2cm}
\def \margin {8} % margin in angles, depends on the radius

\foreach \s in {1,...,\n}
{
  \node[draw, circle] at ({360/\n * (\s - 1)}:\radius) {$\alph{\Labels(\s)}$};
  \draw[->, >=latex] ({360/\n * (\s - 1)+\margin}:\radius) 
    arc ({360/\n * (\s - 1)+\margin}:{360/\n * (\s)-\margin}:\radius);
}
\end{tikzpicture}}

\begin{figure}[htb]
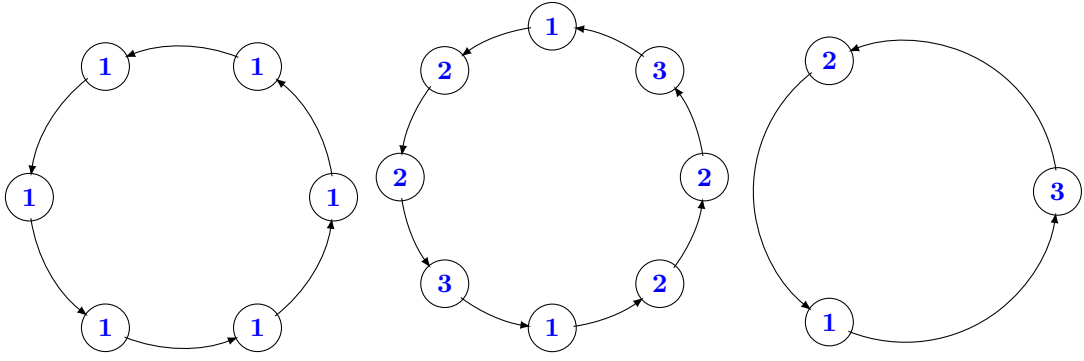

\begin{center}
\necklace{1&1&1&1&1&1}{6}
\necklace{2&3&1&2&2&3&1&2}{8}
\necklace{3&2&1}{3}
\caption{\label{fig:necks}Three necklaces}
\end{center}
\end{figure}

Consider ordered repeating patterns with beads in $\Sigma$ like those in Figure~\ref{fig:necks}. They are known as necklaces. We represent the pattern with the word of the letters forming the repeating cycle which follows the order of the arrows and comes first in alphabetical order. For the three necklaces shown, the repeating pattern is $\alph{1}$, $\alph{1223}$ and $\alph{132}$ respectively.

Not all words represent a repeating pattern in this setup. Words which do are called \emph{Lyndon words}. For example, $\alph{123}$ is a Lyndon word but $\alph{21}$ is not. The empty word is not a Lyndon word. Every single-letter word is a Lyndon word. A word which is comprised of multiple copies of a smaller word is not a Lyndon word. A multi-letter word is a Lyndon word if and only if it comes before all its proper suffixes in alphabetical order. %Also, a multi-letter word is a Lyndon word iff it comes before all its rotations in alphabetical order.
If $w=uv$ is a multi-letter Lyndon word and $v$ is the longest proper suffix of $w$ which is also a Lyndon word, then $u$ is a Lyndon word. 
A formula for the number of Lyndon words of length $m$ on $d$ letters is given by \cite{wikinecklace}.

\section{Bracketed expressions}

\renewcommand{\c}{\mathbin{\raisebox{0ex}{,}}}%\mathbin{\raisebox{0.5ex}{,}}}
\def\myexample{$\left[\,\alph{1}\c\left[\,\left[\,\alph{2}\c\alph{1}\,\right]\c\alph{4}\,\right]\,\right]$}
A bracketed expression on $\Sigma$ is an expression consisting of one or more letters combined with a binary $[\,\cdot\c\cdot\,]$ operator called the Lie bracket. For example \myexample\ or $\alph{3}$ or $[\,\alph{2}\c\alph{1}\,]$.
% The depth of a bracketed expression is defined recursively as 1 for a letter and 1 plus the maximum of the depths of $u$ and $v$ for $[u,v]$. For instance, the depth of \myexample\ is 3. %USELESS CONCEPT
The \emph{foliage} of a bracketed expression is the word formed from its letters, ignoring its brackets and commas. 
The \emph{depth} of a bracketed expression is the length of its foliage, which is equivalent to one plus the number of sets of brackets it contains. For instance, the foliage of \myexample\ is $\alph{1214}$ and its depth is 4. Note that the depth is not in general equal to one plus the maximum number of brackets surrounding a letter, for example consider $\left[\,\left[\,\alph{2}\c\alph{1}\,\right]\c\left[\,\alph{1}\c\alph{4}\,\right]\,\right]$.

%consider using commutator from the physics package?

%Just as we considered the tensor algebra from words above, we can consider the \emph{free algebra} on $\Sigma$ as the set of functions from bracketed expressions to $\mathbb{R}$ or the set of possibly infinite formal sums of linear multiples of bracketed expressions.

On the real vector space of linear combinations of bracketed expressions, we can extend the bracket operation to a bilinear operation. For example \[ \Big[\,\big(3\alph{1}+9\,[\,\alph{3}\c\alph{2}\,]\,\big)\c4\alph{2}\Big]=12\,[\,\alph{1}\c\alph{2}\,]+36\,[\,[\,\alph{3}\c\alph{2}\,]\c\alph{2}\,].\]

The \emph{free Lie%\footnote{pronounced "lee"}
~algebra} is this space but where we consider an element unchanged if it undergoes either of the following transformations. %(It is the set of equivalence classes under the equivalence relation generated by them.)
%\begin{itemize}
%\item $[u\c v] \mapsto -1[v\c u]$
%\item $[u\c [v\c w]] \mapsto [[u\c v]\c w] + [[w\c u]\c v]$
%//\end{itemize}
\[[u\c v] \mapsto -1[v\c u]\qquad \qquad[u\c [v\c w]] \mapsto [[u\c v]\c w] + [[w\c u]\c v]\]

This means that for any elements $u, v, w$ of the free Lie algebra, denoting the zero element by $0$, $[u\c u]=0$, and $u$, $v$ and $w$ obey the \emph{Jacobi identity}: 
\begin{equation}[[u\c v]\c w] + [[w\c u]\c v]+ [[v\c w]\c u]=0.
\end{equation}

For example, as elements of the free Lie algebra, the following two expressions are equal \begin{align*}
-13\,[\alph{1}\c\alph{2}]+&19\,[\alph{3}\c[[\alph{1}\c\alph{2}]\c\alph{4}]\\
13\,[\alph{2}\c\alph{1}]+&19\,[[\alph{3}\c[\alph{1}\c\alph{2}]]\c\alph{4}]+19\,[[\alph{4}\c\alph{3}]\c[\alph{1}\c\alph{2}.]]
\end{align*}

The addition and the multiplication by a scalar are well-defined on the equivalence classes. The free Lie algebra is a real vector space and the Lie bracket is still bilinear.

Level $m$ of the free Lie algebra is the subspace spanned by bracketed expressions of depth $m$. (This is well-defined as the transformations do not change the depth of an expression.) The free Lie algebra is the direct sum of its levels. We define the free Lie algebra up to level $m$, $F^{\underline m}(\Sigma)$, as the direct sum of the first $m$ levels. On $F^{\underline m}(\Sigma)$, we define the Lie bracket operation as the projection of the free Lie algebra's Lie bracket onto $F^{\underline m}(\Sigma)$.

We can define a linear map $\rho$ from the free Lie algebra to $T(\mathbb{R}^d)$ which takes letters to letters and $[u,v]$ to $\rho(u)\rho(v)-\rho(v)\rho(u)$. This definition is consistent, i.e.~well defined on equivalence classes.  For example $\rho([\alph{1}\c[\alph{1}\c\alph{2}]])=\alph{112}+\alph{211}-2\alph{121}$. An element of the image of $\rho$ is called a \emph{lie element} and it never has a component in $\epsilon$. For each $m$, $\rho$ is an injection from $F^{\underline m}(\Sigma)$ to $T^{\underline m}(\mathbb{R}^d)$. We identify an element $u$ of the free Lie algebra with $\rho(u)$.

\subsection{Lyndon basis}
Having a basis of each level of the free Lie algebra is useful, for instance
because it gives a canonical representation of each element. %because it makes it possible to represent each element unambiguously. 
Bracketed expressions of depth $m$ span level $m$ of the free Lie algebra, so it is reasonable to expect a basis in terms of bracketed expressions. Famously, Lyndon words can be used to do this.\footnote{I found \cite{LR} to be the easiest introduction to this.}

Define a function $\sigma$ from Lyndon words to bracketed expressions recursively as follows:
\begin{itemize}
\item If $a\in\Sigma$, $\sigma(a)=a$
\item If $w=uv$ is a multi-letter Lyndon word and $v$ is the longest proper suffix of $w$ which is a Lyndon word, then $\sigma(w)=[\sigma(u)\c\sigma(v)]$
\end{itemize}
Then in fact the image under $\sigma$ of the set of Lyndon words of length $m$ forms a basis of level $m$ of the free Lie algebra. %Aside from the $\sigma$ symbol and consider $F^{\underline m}(\Sigma)$ to be the vector space with basis given by the Lyndon words of length $m$ or less. 
For example, Table~\ref{tab:lynd} shows the elements of the Lyndon basis of $F^{\underline 4}(\{\alph{1},\alph{2}\})$.
\begin{table}[ht]
\begin{center}
\begin{tabular}{r|c|c}
level&$u$&$\sigma(u)$\\
\hline
1&$\alph{1}$&$\alph{1}$\\
1&$\alph{2}$&$\alph{2}$\\
2&$\alph{12}$&$[\alph{1}\c\alph{2}]$\\
3&$\alph{112}$&$[\alph{1}\c[\alph{1}\c\alph{2}]]$\\
3&$\alph{122}$&$[[\alph{1}\c\alph{2}]\c\alph{2}]$\\
4&$\alph{1112}$&$[\alph{1}\c[\alph{1}\c[\alph{1}\c\alph{2}]]]$\\
4&$\alph{1122}$&$[[\alph{1}\c[\alph{1}\c\alph{2}]]\c\alph{2}]$\\
4&$\alph{1222}$&$[[[\alph{1}\c\alph{2}]\c\alph{2}]\c\alph{2}]$\\
\end{tabular}
\caption{\label{tab:lynd}The Lyndon words on two letters up to level 4, and their image under $\sigma$}
\end{center}
\end{table}

If $u$ and $v$ are Lyndon words, then $[\,\sigma(u)\c\sigma(v)\,]$ can be calculated in the Lyndon basis according to the following recursive procedure. (Use only the first relevant step at each stage.)
\begin{itemize}
\item If $u=v$, then $[\,\sigma(u)\c\sigma(v)\,]:=0$
\item If $v\lex u$, then $[\,\sigma(u)\c\sigma(v)\,]:=-\left[\,\sigma(u)\c\sigma(v)\,\right]$
\item If $|u|=1$, then $[\,\sigma(u)\c\sigma(v)\,]:=\sigma(uv)$
\item Let $u=xy$, with $y$ being the longest proper suffix of $u$ which is a Lyndon word. If $v\lex y$, then $[\,\sigma(u)\c\sigma(v)\,]:=[\,\sigma(y)\c[\,\sigma(v)\c\sigma(x)\,]\,]+[\,\sigma(x)\c[\,\sigma(y)\c\sigma(v)\,]\,]$
\item $[\,\sigma(u)\c\sigma(v)\,]:=\sigma(uv)$
\end{itemize}

For example, $[\,\alph{2}\c\sigma(\alph{13})\,]$ is $-\sigma(\alph{132})$, and $[\,\alph{3}\c\sigma(\alph{12})\,]$ is $\sigma(\alph{132})-\sigma(\alph{123})$. Having this procedure, we can do arithmetic between elements of the free Lie algebra cleanly storing all elements in terms of the Lyndon basis. Knowing the longest Lyndon proper suffix of a Lyndon word is seen to be important, and any computer program doing these calculations will be designed with that in mind.

%use a symbol for concatenation product ?
%Use * to identify inessential sections
\section{The Baker-Campbell-Hausdorff formula}
%and so BCH stays in Im(rho)
\def\x{\alph{1}}
\def\y{\alph{2}}
With $\log$ and $\exp$ defined according to the series above, we consider the expression $\log(\exp(\x)\exp(\y))$ expanded up to a fixed level $m$. It is in fact a value in $F^{\underline m}(\left\{\alph1,\alph2\right\})$, in other words it is in the image of $\rho$. For example, up to level 2,
\begin{align*}
\gdef\fr#1#2{\tfrac{#1}{#2}}
\log(\exp(\x)\exp(\y))&=\log\big((1\epsilon+\x+\fr12\x\x)(1\epsilon+\y+\fr12\y\y)\big)
\\&=\log(1\epsilon+\x+\y+\fr12\x\x+\x\y+\fr12\y\y)
\\&=(\x+\y+\fr12\x\x+\x\y+\fr12\y\y)-\fr12(\x\x+\y\y+\x\y+\y\x)
\\&=\x+\y+\fr12\x\y-\fr12\y\x
\\&=\x+\y+\fr12(\x\y-\y\x)\\
&=\x+\y+\fr12[\x\c\y]
\end{align*}
Similarly up to level 3
\begin{align*}
\log(\exp(\x)\exp(\y))&=\log\big((1\epsilon+\x+\fr12\x\x+\fr16\x\x\x)(1\epsilon+\y+\fr12\y\y+\fr16\y\y\y)\big)
\\&=\log(1\epsilon+\x+\y+\fr12\x\x+\x\y+\fr12\y\y+\fr16\x\x\x+\fr12\x\y\y+\fr12\x\x\y+\fr16\y\y\y)
\\&=(\x+\y+\fr12\x\x+\x\y+\fr12\y\y+\fr16\x\x\x+\fr12\x\y\y+\fr12\x\x\y+\fr16\y\y\y)
\\&\qquad-\fr12(\x\x+\y\y+\x\y+\y\x+\x\x\x+\fr32\x\x\y
\\&\qquad\qquad\qquad+\x\y\x+\fr32\x\y\y+\fr12\y\x\x+\y\x\y+\fr12\y\y\x+\y\y\y)
\\&\qquad+\fr13(\x\x\x+\x\x\y+\x\y\x+\x\y\y+\y\x\x+\y\x\y+\y\y\x+\y\y\y)
\\&=\x+\y+\fr12\x\y-\fr12\y\x+\fr1{12}\x\x\y-\fr16\x\y\x+\fr1{12}\x\y\y+\fr1{12}\y\y\x-\fr16\y\x\y+\fr1{12}\y\y\x
\\&=\x+\y+\fr12(\x\y-\y\x)
\\&\qquad+\fr1{12}(\x(\x\y-\y\x)-(\x\y-\y\x)\x)
						  -\fr1{12}((\y\x-\x\y)\y-\y(\y\x-\x\y))
\\&=\x+\y+\fr12[\x\c\y]+\fr1{12}[\x\c[\x\c\y]]-\fr1{12}[[\y\c\x]\c\y]				  
\end{align*}
This calculation always works, although the algebraic manipulation gets much more tedious as $m$ increases.
The expression up to an arbitrary number of terms is called the Baker-Campbell-Hausdorff (BCH) formula. It can be written in any basis for  $F^{\underline m}(\{\alph{1},\alph{2}\})$.
The coefficients can be calculated in the Lyndon basis using the method of \cite{bch}. The authors have saved the results up to level 20 a file called \verb|bchLyndon20.dat| in \cite{bchinfo}, and in practice we just use their values. If $x$ and $y$ are arbitrary in $F^{\underline m}(\Sigma)$, then they can be substituted in as $\x$ and $\y$ to give an expression for $\log(\exp(x)\exp(y))$ as a member of $F^{\underline m}(\Sigma)$. 
For example, in $F^{\underline 4}(\{\alph{1},\alph{2},\alph3\})$,
\begin{align*}
\log(\exp([\alph1,\alph2])\exp(\alph3))&=[\alph1,\alph2]+\alph3+\tfrac12[[\alph1,\alph2],\alph3]-\tfrac1{12}[[\alph3,[\alph1,\alph2]],\alph3].
\end{align*}
Here I could ignore the fourth level of the BCH formula because, seeing as all its terms will have a $\alph1$, I will get terms of depth 5 after the substitution.

\newpage
\section{Signature}
\label{sec:sig} 

A path in $\mathbb{R}^d$ %(with coordinate axes $x_1,x_2,\ldots,x_d$) 
can be described by a continuous map $\gamma:[a,b]\to \mathbb{R}^d$ with $\gamma(t)=(\gamma_{\alph1}(t),\gamma_{\alph2}(t),\ldots,\gamma_{\alph d}(t))$. Its \emph{signature} is an element of $T(\mathbb{R}^d)$, denoted by $X^\gamma_{a,b}$.
The component the signature in each word (its image on each word) is defined inductively on the length of the word. The signature of the empty word is defined as 1. If $w$ is a word and $i\in\{\alph{1},\alph{2},\ldots,\alph{d}\}$ then $X^\gamma_{a,b}(wi)$ is defined as $\int_a^tX^\gamma_{a,t}\gamma_i'(t)\,dt$. 
The  restriction of the signature to words of length $m$ is called the $m$th \textbf{level} of the signature. It contains $d^m$ values. Written in Riemann-Stieltjes form, they are the values of integrals of the form 
\begin{equation}
\int_a^b\int_a^{t_1}\dots\int_a^{t_{m-2}}\int_a^{t_{m-1}}\,d\gamma_{i_1}(t_m)\,d\gamma_{i_2}(t_{m-1})\,\dots\,d\gamma_{i_{m-1}}(t_2)\,d\gamma_{i_m}(t_1),
\end{equation}
 where each $i_j$ is allowed to range over values in $\{\alph1,\alph2,\ldots,\alph d\}$. For instance, level 1 of the signature is the total displacement of the path. If $d\ge3$ then the $\alph{231}$ component of the signature of is the value of the integral \begin{align}
\int_a^b\int_a^{t}\int_a^{s}\gamma'_{\alph2}(r)\,dr\, \gamma'_\alph3(s)\,ds\,\gamma'_{\alph1}(t)\,dt &=\iiint_{a<r<s<t<b}\gamma'_{\alph2}(r)\gamma'_\alph3(s)\gamma'_{\alph1}(t)\,dV\\&=\int_a^b\int_a^{t}\int_a^{s}\,d\gamma_{\alph2}(r)\,d\gamma_{\alph3}(s)\,d\gamma_{\alph1}(t)
\\&=\int_a^b\int_r^{b}\int_s^{b}\,d\gamma_{\alph1}(t)\,d\gamma_{\alph3}(s)\,d\gamma_{\alph2}(r).
\end{align}
We are interested in the signature up to level $m$, $X^{\gamma,\underline m}_{a,b}$, which is an element of $T^{\underline m}(\mathbb{R}^d)$.

%The signature of a straight line is symmetric in 
If $a<b<c$ then the result (from \cite{chen}) known as \textbf{Chen's identity} states that \begin{equation}
X^\gamma_{a,c}(i_1i_2\ldots i_n) %$
=\sum_{j=0}^nX^\gamma_{a,b}(i_1i_2\ldots i_{j})X^\gamma_{b,c}(i_{j+1}i_{j+2}\ldots i_n).
\end{equation}
(The products indicated with ellipses can be empty, indicating the empty word, on which any signature takes the value 1.) For example
\\
\begin{minipage}{0.8\textwidth}
\begin{align*}
X^\gamma_{a,c}(\alph1\alph2)&=\int_a^c\int_a^s\gamma'_\alph1(r)\,dr\,\gamma'_\alph2(s)\,ds 
\\&=\left[\int_b^c\int_b^s+\int_b^c\int_a^b+\int_a^b\int_a^s\right]\gamma'_\alph1(r)\,dr\,\gamma'_\alph2(s)\,ds 
\\&=X^\gamma_{b,c}(\alph1\alph2)+X^\gamma_{a,b}(\alph1)X^\gamma_{b,c}(\alph2) + X^\gamma_{a,b}(\alph1\alph2).
\end{align*}
\end{minipage}%
\begin{minipage}{0.2\textwidth}
\begin{tikzpicture}[scale=0.8]
\draw[->] (0,0) -- (2.5,0) node[anchor=west] {$s$};
\draw[->] (0,0) -- (0,2.5) node[anchor=south] {$r$};
\draw (0,0) -- (0,-2pt) node[anchor=north] {$\vphantom{b}a$};
\draw (1,0) -- (1,-2pt) node[anchor=north] {$b$};
\draw (2,0) -- (2,-2pt) node[anchor=north] {$\vphantom{b}c$};
\draw (0,0) -- (-2pt,0) node[anchor=east] {$a$};
\draw (0,1) -- (-2pt,1) node[anchor=east] {$b$};
\draw (0,2) -- (-2pt,2) node[anchor=east] {$c$};
\fill[black!80] (0,0) -- (1,0) -- (1,1);
\fill[black!50] (1,1) -- (1,0) -- (2,0) -- (2,1);
\fill[black!30] (1,1) -- (2,1) -- (2,2);
\end{tikzpicture}
\end{minipage}

Restricting up to level $m$, this means that 
\begin{equation}\label{eq:chen}
X^{\gamma,\underline m}_{a,c}=X^{\gamma,\underline m}_{a,b}X^{\gamma,\underline m}_{b,c},
\end{equation}
using the concatenation product in $T^{\underline m}(\mathbb{R}^d)$.

Also, if $\gamma$ is a straight line then \begin{equation}X^\gamma_{a,b}(i_1i_2\ldots i_n)=\frac1{n!}\prod_{j=1}^n(\gamma_{i_j}(b)-\gamma_{i_j}(a)).\end{equation} 
Restricting up to level $m$, this means that, for the straight path, 
\begin{align}\label{eq:line}
X^{\gamma,\underline m}_{a,b}&=\exp\Big((\gamma_{\alph1}(b)-\gamma_{\alph1}(a))\alph{1} +\dots+ (\gamma_{\alph d}(b)-\gamma_{\alph d}(a))\alph{d}\Big)
=\exp\big(\gamma(b)-\gamma(a)\big)
\end{align}
Equations \ref{eq:chen} and \ref{eq:line} between them make it easy to calculate elements of the signature of a piecewise linear path up to any given level. The code required to do this is simple, although it requires more operations than direct log signature calculations.

\section{Log signature}
\label{sec:logsig}
Let $S$ be the set which consists of the signature up to level $m$ of every path in $\mathbb{R}^d$.
$S$ is not the whole of the vector space $T^{\underline m}(\mathbb{R}^d)$, but rather a lower dimensional manifold. $S$ is known as the \emph{group-like elements} of $T^{\underline m}(\mathbb{R}^d)$.
In fact, $S$ is the image under $\exp$ of the image of $\rho$. The \emph{log signature} of $\gamma$ up to level $m$ is $Y^{\gamma,\underline m}_{a,b}=\rho^{-1}(\log (X^{\gamma,\underline m}_{a,b}))$, an element of $F^{\underline m}(\{\alph{1},\dots,\alph{d}\})$. When the term `log signature' is used to describe a fixed array of numbers, of course we mean its components in some basis, for example the Lyndon basis.

For speed, it may be best to calculate the log signature explicitly without leaving the free Lie algebra. Here are the two rules which enable this to happen. Equation~\ref{eq:line} tells us that if $\gamma$ is a straight line then its log signature is independent of $m$: \begin{equation}Y^{\gamma,\underline m}_{a,b}=\gamma(b)-\gamma(a).\end{equation}
From Equation~\ref{eq:chen} it follows that for  $a<b<c$ 
\begin{equation}\label{eq:chenbch}
Y^{\gamma,\underline m}_{a,c}=\log\Big(\exp\big(Y^{\gamma,\underline m}_{a,b}\big)\exp\big(Y^{\gamma,\underline m}_{b,c}\big)\Big).
\end{equation}
Thus the log signature of the concatenation of two paths can be calculated from the paths using the Baker-Campbell-Hausdorff formula, with the $\x$ and $\y$ there above replaced with the log signatures of the two paths.

\section{Sizes}
The size of the signature and log signature up to level $m$ for various values of $d$ and $m$ is as given in Table~\ref{tab:sizes}. The fixed 1 coefficient of $\epsilon$ in the signature is omitted as it does not carry information.
\begin{table}[ht]
\begin{center}
\def\l#1{&\textbf{#1}}
\gdef\m#1{&\multicolumn{2}{c|}{$d=#1$}}
\begin{tabular}{|c|rl|rl|rl|rl|rl|}
\hline
\m 1 \m 2 \m 3 \m 4 \m 5\\
\hline
$m=1$&1\l1&2\l2&3\l3&4\l4&5\l5\\
$m=2$&2\l1&6\l3&12\l6&20\l{10}&30\l{15}\\
$m=3$&3\l1&14\l5&39\l{14}&84\l{30}&155\l{55}\\
$m=4$&4\l1&30\l8&120\l{32}&340\l{90}&780\l{205}\\
$m=5$&5\l1&62\l{14}&363\l{80}&1,364\l{294}&3905\l{829}\\
$m=6$&6\l1&126\l{23}&1,092\l{196}&5,460\l{964}&19,530\l{3,409}\\
$m=7$&7\l1&254\l{41}&3,279\l{508}&21,844\l{3,304}&97,655\l{14,569}\\
$m=8$&8\l1&510\l{71}&9,840\l{1318}&87,380\l{11,464}&488,280\l{63,319}\\
$m=9$&9\l1&1022\l{127}&29,523\l{3,502}&349,524\l{40,584}&244,1405\l{280,319}\\
$m=10$&10\l1&2046\l{226}&88,572\l{9,382}&\small 1,398,100\l{145,338}&\small 12,207,030\l{\small 1,256,567}\\
$m=11$&11\l1&4094\l{412}&265,719\l{25,486}&\small 5,592,404\l{526,638}&\small 61,035,155\l{\small 5,695,487}\\
$m=12$&12\l1&8190\l{747}&797,160\l{69,706}&\small $\,$\llap{2}2,369,620\l{\small 1,924,37\rlap{8}\,}&\small $\,$\llap{3}05,175,780\l{\small 26,039,18\rlap{7}\,}\\
\hline
\end{tabular}
\caption{\label{tab:sizes}The sizes of signatures (plain) and log signatures (bold)}
\end{center}
\end{table}

The size of the signature is $\frac{d(d^m-1)}{d-1}$. A formula for the size of each level of the log signature is given in \cite{wikinecklace}.

\section{The code}
\begin{center}
	\begin{minipage}{0.7\textwidth}
		\noindent\textbf{Note:} This section refers to illustrative code which can be found at \url{https://github.com/bottler/LogSignatureDemo}. This document predates the Python package \verb|iisignature| which provides ready access to these calculations.
	\end{minipage}
\end{center}
\label{sec:code}
If the whole log signature up to level $m$ is being used as a representation in a machine learning algorithm, then often the choice of basis should not matter. Typically the aim is to have fast code to calculate the log signature up to a fixed level in the same basis for many different paths of a fixed dimension. The log signature of a single segment of a path is easy to calculate. The python script is used to write the code which constructs signatures of multi-segment paths.

The python script represents basis elements of the free Lie algebra by objects of class \verb|LyndonWord|, and elements of the free Lie algebra by objects of type \verb|Polynomial|, which holds a dictionary from \verb|LyndonWord| to \emph{numbers}. By using an algebra of polynomials in arbitrary strings as the \emph{numbers} (classes \verb|N_constantNum| and \verb|N_polynomial|), it can calculate expressions for the components of the log signature of the concatenation of two paths given their log signatures. The script writes functions to do this in the code which it outputs. When the second path is a line segment, the expressions are considerably simpler, and a function for this can also be written.
The user needs to pay attention to areas in the code where the comment \verb|USER CONTROL| occurs, to determine such things as which $m$ and $d$ to consider and where the data from \cite{bchinfo} has been downloaded.
The code produced is not affected by which basis the BCH formula is expressed in, as all the BCH's brackets get expanded.
\subsection*{Using the generated code}
The generated code is used in \verb|C++| to calculate the log signature of a given path. The generated file \verb|bch.h| will need to be \verb|#include|d, and the \verb|bch.cpp| will need to be built in the client project. Note that \verb|bch.cpp| has no other dependencies and if high levels are included will take a while to compile, so it is recommended to design the client's build procedure around not recompiling it every build. Here is a quick example of how to use it: If \verb|p| is a \verb|std::vector<std::array<float,3>>| containing the coordinates visited by a 3d path, then the following code will make \verb|logsig| be a \verb|std::array| containing its log signature up to level 6.
\begin{quote}
\begin{verbatim}
LogSignature<3,6> logsig{};
if(p.size()>1)
  for(int j=0; j<3; ++j)
    logsig[j]=p[1][j]-p[0][j];
for(size_t i=2; i<p.size(); ++i){
  Segment<3> displacement;
  for(int j=0; j<3; ++j)
    displacement[j]=p[i][j]-p[i-1][j];
  joinSegmentToSignatureInPlace(logsig,displacement);
}
\end{verbatim}
\end{quote}

\section{Other bases}
\label{sec:other}
This open source \verb|C++| library \cite{coropa} is able to work with elements of the free Lie algebra, which are represented as variables of type \verb|alg::lie|. Like we do, the library stores them in terms of a basis. The calculation of the basis is in the file \verb|libalgebra/lie_basis.h|, and examples of the basis which it calculates can be generated using the form on the website. The basis used is not the Lyndon basis described earlier.

In \cite{hall1950} is given a general method of constructing a basis of the free Lie algebra once an order on bracketed expressions with a given depth has been specified. (This is extended to an order on all bracketed expressions by saying that expressions with higher depth are greater than those with lower depth). This is called a (generalized) Hall basis. If the order is alphabetical on the foliage, then the basis obtained is the Lyndon basis.

Coropa uses a different order. For Coropa, if two different bracketed expressions $[u\c v]$ and $[x\c y]$ have the same length then $[u\c v]<[x\c y]$ if $u<x$ or both $u=x$ and $v<y$. In fact, Coropa uses a convention in which all brackets are reversed. The effect is that terms in even levels (which have an odd number of brackets) are negated. The basis elements for 2 dimensions begin $\alph1$, $\alph2$, $[\alph{1}\c\alph{2}]$,
$[\alph{1}\c[\alph{1}\c\alph{2}]]$,
$[\alph{2}\c[\alph{1}\c\alph{2}]]$.

If the variable \verb|match_coropa| is set to \verb|True| in the python code, then it will use this basis. Note that in this case the class \verb|LyndonWord| represents a basis element, which is not actually a Lyndon word. The generated code is a bit longer when this is done; it is not entirely clear why this is.
 
In \cite{bch} an algorithm is given for expressing the BCH formula using any generalized Hall basis with $d=2$. One of the example bases which is given is called \emph{the (classical) Hall basis}. This is the same as Coropa's basis on 2 words, except that the brackets are not reversed. This means that basis elements of odd depth are the same, and those of even depth are negated. The basis elements for 2 dimensions begin $\alph1$, $\alph2$, $[\alph{2}\c\alph{1}]$,
$[[\alph{2}\c\alph{1}]\c\alph{1}]$,
$[[\alph{2}\c\alph{1}]\c\alph{2}]$.
\section{A visualisation}
\label{sec:vis}
It is natural to want to visualise the log signature, and here is a small idea in that direction.
In the specific case $d=2$, $m=4$, the log signature has eight components, and a two dimensional path with four segments has eight degrees of freedom. The file \verb|view.m| provides an interactive Mathematica~10 visualisation of the relationship between a path made of four segments and its log signature. It depends on another file \verb|bch.m|, which has been generated by additional functionality in the python script, which defines a single function returning the log signature of a path defined by four displacements. The visualisation should appear when \verb|view.m| is run. The log signature appears as widgets on the right of the graph of the path. Because level 1 and level 3 of the log signature are two dimensional, they are represented by 2d controls. The other log signature components are controlled separately. When `solve' is ticked, you can gently move the widgets to change the components of the log signature and see the path move. Note that the 12 component controls the signed area enclosed by the path. The calculation of solutions is not perfect, but is enough to get a general picture. When `solve' is unticked, you can drag the locators to change the path, and see the corresponding log signature elements move. Both files can be found at 
\url{https://github.com/bottler/iisignature/tree/master/examples/Mathematica}.% is where both files can be found.

\vspace{-1.4em}%so we don't have a hanging couple of bibliography items

\printbibliography
\end{document}